\documentclass[11pt,oneside]{amsart}
\usepackage{amsmath,amssymb,graphicx}
\newtheorem{thm}{Theorem}[section]
\newtheorem{lem}[thm]{Lemma}
\newtheorem{cor}[thm]{Corollary} 

\theoremstyle{definition}
\newtheorem{defn}[thm]{Definition}
\newtheorem{rem}[thm]{Remark}
\newtheorem{term}[thm]{Terminology}
\newtheorem{exmp}[thm]{Example}

\newcommand{\blackboard}[1]{\ensuremath{\mathbb{#1}}}
\newcommand{\script}[1]{\ensuremath{\mathcal{#1}}}
\newcommand{\smallcaps}[1]{\textrm{\textsc{#1}}}
\newcommand{\sph}{\blackboard{S}}
\newcommand{\R}{\blackboard{R}}

\newcommand{\cat}{\smallcaps{CAT}}

\newcommand{\bt}{\begin{tabular}}
\newcommand{\et}{\end{tabular}}

\begin{document}

\title[Word-Hyperbolic $3$-manifolds]
  {Combinatorial conditions that imply\\ 
word-hyperbolicity for $3$-manifolds}

\author[M.~Elder]{Murray Elder$\!{ }^1$}
\address{Dept. of Mathematics\\
	Tufts University\\
	Medford, MA 02155}
\email{murray.elder@math.tamu.edu}

\author[J.~McCammond]{Jon McCammond$\!{ }^1$}
\address{Dept. of Mathematics\\
	U. C. Santa Barbara\\
	Santa Barbara, CA 93106}
\email{jon.mccammond@math.ucsb.edu}

\author[J.Meier]{John Meier}
\address{Dept. of Mathematics\\
	Lafayette College\\
	Easton, PA 18042}
\email{meierj@lafayette.edu}

\date{\today}

\begin{abstract} 
Thurston conjectured that a closed triangulated $3$-manifold in which
every edge has degree $5$ or $6$, and no two edges of degree $5$ lie
in a common $2$-cell, has word-hyperbolic fundamental group.  We
establish Thurston's conjecture by proving that such a manifold admits
a piecewise Euclidean metric of non-positive curvature and the
universal cover contains no isometrically embedded flat planes.  The
proof involves a mixture of computer computation and techniques from
small cancellation theory.
\end{abstract}

\footnotetext[1]{Partially supported under NSF grant DMS-0101506}

\maketitle


\section{Introduction }
In this article we show that a class of closed triangulated
$3$-manifolds can be assigned a metric of non-positive curvature.  In
addition to proving a conjecture of Thurston, our main result
illustrates the way in which the computer program developed by the
first and second authors (\cite{ElMc-3d}) can be used in conjunction
with combinatorial methods to establish non-trivial results about
$3$-manifolds.  The class of triangulations we consider are defined as
follows.

\begin{defn}[$5/6^*$-triangulations] 
Let $M$ be a closed triangulated $3$-manifold and recall that the
\emph{degree} of an edge is the number of closed tetrahedra which
contain it.  If every edge in $M$ has degree $5$ or $6$ then this is a
\emph{$5/6$-triangulation} of $M$.  The triangulation is called a
\emph{$5/6^*$-triangulation} if each $2$-cell in $M$ contains at most
one edge of degree $5$.
\end{defn}

Thurston conjectured that every closed $3$-manifold that admits a
$5/6^*$-triangulation has word-hyperbolic fundamental group.  We prove
a slightly stronger version of this conjecture.

\begin{thm}[Main Theorem]\label{thm:main} 
Every $5/6^*$-triangulation of a closed $3$-manifold $M$ admits a
piecewise Euclidean metric of non-positive curvature, where the
universal cover $\widetilde{M}$ contains no isometrically embedded
flat planes.  As a consequence, $\pi_1(M)$ is word-hyperbolic.
\end{thm}

Such triangulations are not as special as they might appear.  Cooper
and Thurston show that every closed $3$-manifold can be cellulated so
that each $3$-cell is a cube and each edge degree is $3$, $4$, or $5$
\cite{CoTh88}.  Using similar techniques, Noel Brady and the last two
authors show that every closed $3$-manifold admits a triangulation
where each edge degree is $4$, $5$, or $6$ \cite{BrMcMe-degrees}.  The
dual concept --- $5/6$-triangulations where $2$-cells contain at most
one edge of degree $6$ --- is called a \emph{foam} and is of interest
in chemistry (see \cite{Su}).

\medskip

\noindent 
{\bf Structure of the paper:} In Sections~\ref{sec:diagrams}
and~\ref{sec:small} we use ideas resembling small cancellation theory
to investigate the $2$-sphere triangulations that arise as vertex
links in $5/6^*$-triangulations of $3$-manifolds.  In
Section~\ref{sec:algorithm} we review the general algorithm for
determining curvature properties in $3$-dimensional metric polyhedral
complexes and establish an improved version of this algorithm that
uses a lemma of Bowditch to greatly simplify the calculations.  We
describe the piecewise Euclidean metric we assign to a
$5/6^*$-triangulation of a $3$-manifold in Section~\ref{sec:output},
and in Section~\ref{sec:main} we compare the computational results for
this metric with the combinatorial restrictions proved in
Section~\ref{sec:small}.  Their mismatch enables us to establish our
main result.

\section{Diagrams and duals}\label{sec:diagrams}
In this section we recall some standard definitions and introduce
soccer diagrams, which are closely associated with
$5/6^*$-triangulations of $3$-manifolds.  For background on disc
diagrams see \cite{fans}.

\begin{defn}[Disc diagram]
A \emph{disc diagram} $D$ is a contractible combinatorial $2$-complex
together with a specific planar embedding $D\to \R^2$.  If $D$ is
homeomorphic to a disc, $D$ is \emph{non-singular}, otherwise it is
\emph{singular}.
\end{defn}

\begin{defn}[Internal dual]\label{def:dual}
Let $D$ be a disc diagram in which all vertices of degree~$2$ lie on
the boundary cycle of $D$.  The \emph{internal dual} of $D$ is a
subspace of $D$ which consists of a $0$-cell at the center of each
$2$-cell of $D$, a $1$-cell passing through each internal $1$-cell of
$D$ connecting the centers of the $2$-cells on either side, and a
$2$-cell for each interior $0$-cell $v$.  See Figure~\ref{fig:min} for
an illustration.  A similar definition can be given when $D$ is a
triangulation of a $2$-sphere, which results in a \emph{dual
cellulation}.
\end{defn}

The following lemma records the basic properties of internal duals.
See Lemma~5.6 in \cite{fans} for a detailed proof.

\begin{lem}[Dual properties]\label{lem:dual}
If $D$ is a non-singular disc diagram and $E$ is its internal dual,
then $E$ is a contractible, but possibly singular, disc diagram.
Moreover, if $R$ is a $2$-cell in $D$ and $v$ is the corresponding
$0$-cell in $E$, then the number of components of $\partial R \cap
\partial D$ equals the number of components of the link of $v$.
\end{lem}

\begin{term}[Paths, loops, vertex degrees]  Throughout this paper, 
\emph{paths} and \emph{loops} in a cell complex are \underbar{edge}
 paths and loops.  The \emph{degree} of a vertex $v$ is the number of
 edges incident with $v$.
\end{term}

\begin{defn}[Sphere triangulation]
If $M$ is a $5/6^*$-triangulation of a $3$-manifold, the link of a
vertex in $M$ is a triangulation of a $2$-sphere in which each vertex
has degree $5$ or $6$ and no two vertices of degree $5$ are connected
by an edge.  A triangulated $2$-sphere with both of these properties
is called a \emph{$5/6^*$-triangulation} of a $2$-sphere.
\end{defn}

\begin{defn}[Soccer diagram]
Let $D$ be a $5/6^*$-triangulation of a $2$-sphere and let $E$ denote
its dual cellulation.  Since vertices in $D$ have degree $5$ or $6$,
the dual consists of pentagons and hexagons; since $D$ is a
triangulation, every vertex in $E$ has degree $3$; and since no two
vertices degree $5$ in $D$ are connected by an edge, no two pentagons
in $E$ share a common side.  The cell structure $E$ is called a
\emph{soccer tiling} of the $2$-sphere since the standard tiling of a
soccer ball is a simple example with all of these properties.  A
subcomplex of a soccer tiling of a $2$-sphere, which is homeomorphic
to a disc, is called a \emph{soccer diagram}.  Notice that the
embedding of a soccer diagram into the $2$-sphere determines a natural
planar embedding as well.  Thus soccer diagrams are disc diagrams and
the definition further implies that they are non-singular.
\end{defn}  

\begin{defn}[Left/right turn] 
Let $P$ be an immersed directed path in  a soccer
tiling.  If the sphere is oriented then there is a
well-defined notion of a \emph{left/right turn} (the $2$-cells are
thought of as convex and approximately regular).  These
are the only possibilities since every vertex has degree $3$.  The
\emph{turn pattern} for $P$ is the sequence of left and right turns.
For soccer diagrams we define a turn pattern by traversing the
boundary cycle counterclockwise in the induced orientation.
(Since every vertex in $\partial D$ has a connected
link and is of degree $2$ or $3$, the boundary vertices of degree~$2$
correspond to left turns and those of degree~$3$ to right turns.)
Finally, let $n_l$ and $n_r$ denote the total number of left and right
turns, respectively.  The difference $n_l-n_r$ is the
\emph{combinatorial turning angle}.
\end{defn}

The combinatorial turning angle is directly related to the number of
pentagons $D$ contains.

\begin{lem}[Turns and pentagons]\label{lem:pent}
If $D$ is a soccer tiling of a $2$-sphere then $D$ contains exactly
twelve pentagons.  Moreover, if $D$ is a soccer diagram containing
exactly $p$ pentagons, the combinatorial turning angle is $6-p$. 
Hence $p\leq 6$ implies $n_l\geq n_r$.
\end{lem}

\begin{proof} 
We sketch a proof using the Combinatorial Gauss-Bonnet Theorem:
For any angle assignment, the sum of the vertex 
curvatures plus the sum of
the face curvatures is always $2\pi$ times the Euler
characteristic. (See for example 
\cite[Section~4]{fans}.) If we assign an angle of $2\pi/3$ to each corner of
each $2$-cell, the hexagons have curvature $0$, the pentagons have
curvature $\pi/3$, the internal vertices have curvature $0$, the
vertices located at left turns have curvature $\pi/3$ and the vertices
located at right turns have curvature $-\pi/3$.  Thus
$2\pi\cdot 2 = (\pi/3)\cdot p$ for a soccer tiling,
 and $2\pi\cdot 1 =
(\pi/3) \cdot (p + n_l - n_r)$ for a soccer diagram. 
Dividing by
$\pi/3$ and rearranging yields the results.
\end{proof}

\begin{rem}[Limiting pentagons]\label{rem:pent}
If $P$ is a simple closed loop embedded in a soccer tiling $D$, then
it bounds two soccer diagrams $D_1$ and $D_2$ whose intersection is
$P$ and whose union is $D$.  Since $D$ contains only
 twelve pentagons, $P$ always bounds some soccer diagram with at most six
pentagons.
\end{rem}

\begin{defn}[Exposed path]\label{def:exposed}
If $D$ is a soccer diagram whose turn pattern contains $i-1$
consecutive left turns, then $D$ contains a $2$-cell $R$ such that
$\partial R$ and $\partial D$ share a path of length $i$.  We call
this an \emph{exposed path of length $i$} and $R$ is the $2$-cell
which is \emph{exposed}.
\end{defn}

\begin{cor}[Alternating turns]\label{cor:pent}
If $D$ is a soccer diagram with at most six pentagons, then either $D$
contains an exposed path of length $3$, or the turn pattern is
$(rl)^i$ for some $i$ and $D$ contains exactly six pentagons.
\end{cor}

\begin{proof}
By Lemma~\ref{lem:pent}, $n_l\geq n_r$.  If $n_l>n_r$, two left turns
must be adjacent, creating an exposed path of length $3$.  On the other
hand, if $n_l=n_r$, then the only way to avoid adjacent left turns is
for the left and right turns to alternate.  This, in turn, implies
exactly six pentagons by Lemma~\ref{lem:pent}.
\end{proof}

\section{Small diagrams}\label{sec:small}
In this section we prove a key technical result about soccer diagrams
with short boundary cycles, Theorem~\ref{thm:1214}.  We begin by
examining two processes by which soccer diagrams can be decomposed
into smaller soccer diagrams: cut paths and double duals.

\begin{figure}[htbp]
\includegraphics{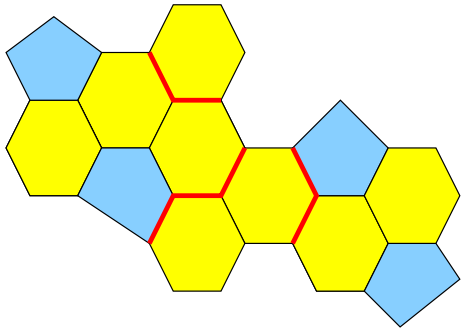}
\caption{A soccer diagram with several cut paths marked by thick
lines.\label{fig:cut}}
\end{figure}

\begin{defn}[Cut path]\label{def:cut}
A soccer diagram $D$ has a {\em cut path} if there are two soccer
diagrams $D_1, D_2$ such that $D=D_1\cup D_2$ and $D_1\cap D_2$ is a
simple path $P$.  See Figure~\ref{fig:cut}.    If $D$ has at
least one cut path, then it has a cut path of minimum length; in
Figure~\ref{fig:cut} this minimum is $2$.
\end{defn}

\begin{lem}[Cut paths exist]\label{lem:cut}
If $D$ is a soccer diagram with at least two $2$-cells then $D$ has a
cut path.  If, in addition, $D$ has at most six pentagons then $D$ has
a cut path of length at most~$4$.
\end{lem}

\begin{proof}
Choose a $2$-cell $R$ such that an edge of $\partial R$ is contained
in $\partial D$.  Since $D$ has more than one $2$-cell, $\partial D$
and $\partial R$ are distinct and there exists a subpath of $\partial
R$ which starts and ends in $\partial D$ and otherwise is contained in
the interior of $D$.  This is a cut path whose length is at most $5$.
If the shortest cut path has length $5$ then every $2$-cell containing
a edge of $\partial D$ must be a hexagon with exactly one edge in
$\partial D$, but this would imply $\partial D$ contains only right
turns, contradicting Lemma~\ref{lem:pent}.
\end{proof}

\begin{figure}[htbp]
\bt{ccc}
\includegraphics[scale=1.4]{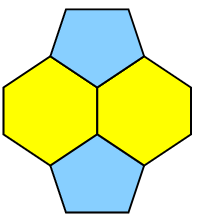}& &
\includegraphics[scale=1.4]{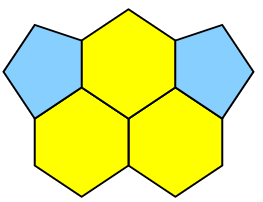}\\
\et
\caption{Soccer diagrams with few $2$-cells and no exposed path of
length $4$.
\label{fig:1214}}
\end{figure}

\begin{lem}[Small diagrams]\label{lem:few}
If $D$ is a soccer diagram with at most five $2$-cells, then $D$
contains a cut path of length at most $2$.  If in addition, $D$ does
not contain an exposed path of length at least $4$, then $D$ is one of
the two disc diagrams shown in Figure~\ref{fig:1214}.
\end{lem}

\begin{proof}
Rather than analyze $D$ directly, it is easier to analyze its internal
dual $E$.  By  Lemma~\ref{lem:dual}, $E$ is a (possibly
singular) disc diagram with at most $5$ vertices and all of
the possibilities for $E$ are shown in Figure~\ref{fig:few}.  Notice
that every possibility for $E$ contains either a vertex of degree~$1$,
or a triangle with a vertex of degree~$2$, and these lead to cut
paths of length~$1$ and~$2$ in $D$.

\begin{figure}[htbp]
\includegraphics[scale=1.4]{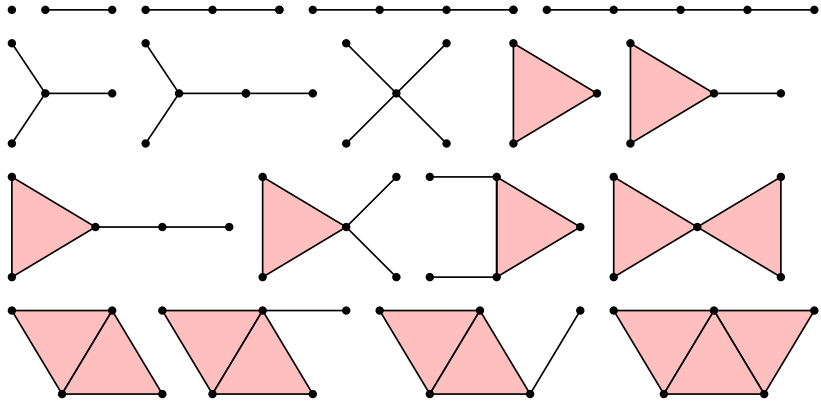}
\caption{Simplicial disc diagrams with at most $5$ vertices.}
\label{fig:few}
\end{figure}

To see the second assertion, notice that a vertex in $E$ of
degree~$1$  corresponds to a $2$-cell $R$ in $D$ such that only
one edge of $\partial R$ lies in the interior of $D$, hence
the remaining edges form an exposed path of length $\ge 4$.
  Similarly, if $E$ contains a triangle in which two of
its vertices have degree~$2$, then these vertices correspond to
$2$-cells in $D$, at least one of which must be a hexagon.  This
hexagon also creates an exposed path of length $4$ in $D$.

There are only two diagrams in Figure~\ref{fig:few} that contain
neither vertices of degree~$1$ nor triangles with two vertices of
degree~$2$ (i.e. the diagrams in the lower left-hand corner and lower
right-hand corner).  In both cases, the vertices in $E$ of degree~$2$
must correspond to pentagons in $D$ in order to avoid exposed paths of
length~$4$, and the remaining vertices must correspond to hexagons
since they share sides with the pentagons.  Thus there are exactly two
possibilities for $D$, and these are the ones shown in
Figure~\ref{fig:1214}.
\end{proof}

The remainder of the section is devoted to showing that soccer
diagrams with more than five $2$-cells and at most six pentagons have
long boundary cycles (Lemma~\ref{lem:six}).  The proof of
Lemma~\ref{lem:six} proceeds by induction and the following lemma
provides the basis step.

\begin{lem}[Minimum length]\label{lem:min}
Let $D$ be a soccer diagram with $k$ $2$-cells and at most six
pentagons.  For $k=1,2,3,4,5,6$, the minimal length of $\partial D$ is
$5,9,11,12,14,15$, respectively.
\end{lem}

\begin{proof}
In Figure~\ref{fig:min} we exhibit soccer diagrams that realize these
values, so the only question is whether there are 
diagrams with smaller boundary
cycles.  Let $E$ be the internal dual of $D$.  If $k\leq
5$, then $E$ is one of the diagrams listed in Figure~\ref{fig:few} and
it is straightforward to enumerate all of the possibilities for $D$
given a specific $E$ and to calculate that the diagrams shown in
Figure~\ref{fig:min}  have minimal length boundary cycles for
these values of $k$. 

\begin{figure}[htbp]
\includegraphics[scale=1.4]{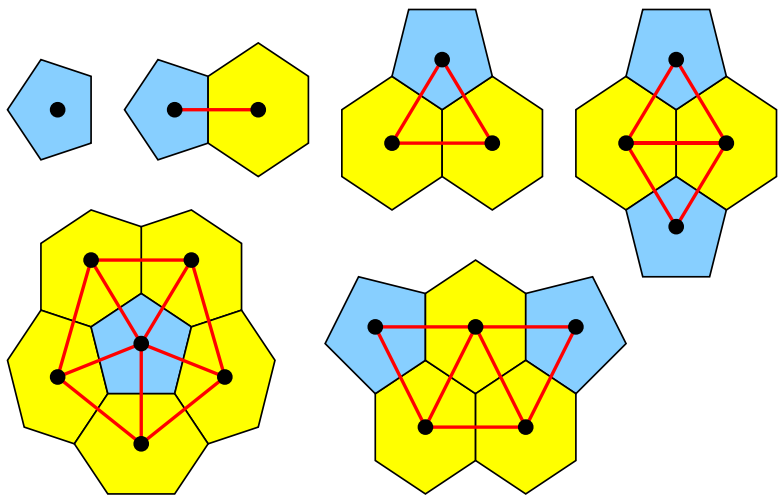}
\caption{Soccer diagrams with $k$ $2$-cells and minimal boundary
lengths for $k$ equals $1$ up to $6$. In each case, the $1$-skeleton
of its internal dual has been superimposed.}
\label{fig:min}
\end{figure}

Now suppose $k=6$.  If $E$ contains either a vertex of degree~$1$ or a
triangle with a vertex of degree~$2$, then $D$ contains a $2$-cell $R$
which is separated from the rest of $D$ by a cut path of length at
most $2$.  Removing $R$ and the exposed portion of $\partial R$ from
$D$ creates a new soccer diagram $D'$ with exactly five $2$-cells, and
$|\partial D'| \geq 14$.  Since $R$ is either a pentagon or a hexagon,
reattaching $R$ to $D'$ shows that $|\partial D| \ge 15$.  The only
six-vertex simplicial disc diagram $E$ that does not contain a vertex
of degree~$1$ or a triangle with a vertex of degree~$2$ is the one
whose $1$-skeleton is superimposed on the soccer diagram in the lower
left-hand corner of Figure~\ref{fig:min}.  For this $E$, the internal
vertex must correspond to a pentagon thereby forcing all of the other
vertices to correspond to hexagons.  In other words, the soccer
diagram shown is the only soccer diagram whose internal dual is $E$,
and the inequality follows.
\end{proof}

The second process we  investigate is the process of taking
double duals.

\begin{figure}[htbp]
\includegraphics[scale=1.2]{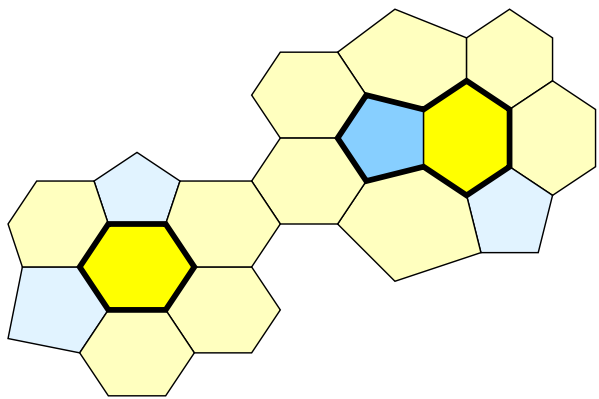}
\caption{A soccer diagram whose double dual is not a soccer diagram
because it is disconnected.\label{fig:double1}}
\end{figure}

\begin{figure}[htbp]
\includegraphics[scale=1.2]{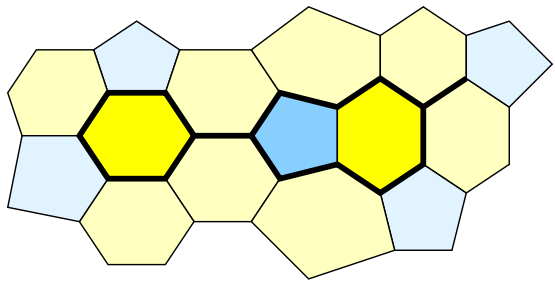}
\caption{A soccer diagram whose double dual is not a soccer diagram
because it is nonsingular.\label{fig:double2}}
\end{figure}

\begin{rem}[Double duals]
Let $D$ be a soccer diagram and let $D'$ be the internal dual of the
internal dual of $D$ (i.e. its \emph{double dual}).  While it is true
that $D'$ is a subcomplex of $D$ (it is essentially $D$ minus the open
star of its boundary), it is not automatically true that $D'$ itself
is a soccer diagram, since $D'$
need not be connected and  it
might be singular.  These situations are illustrated in
Figures~\ref{fig:double1} and \ref{fig:double2}.
\end{rem}

\begin{lem}[Double duals]\label{lem:double}
Let $D$ be a soccer diagram and let $D'$ be its double dual.  If $D'$
is disconnected, then $D$ contains a cut path of length at most $2$.
If $D'$ is connected, but singular, then the turn pattern for $D$
contains either two consecutive left turns or two consecutive right
turns.
\end{lem}

\begin{proof}
Let $E$ denote the dual of $D$.  If $D'$ is disconnected, there are
triangles in $E$ that cannot be connected by a sequence of triangles
so that successive triangles share a common side.  Since $E$ itself is
connected (Lemma~\ref{lem:dual}), this implies that $E$ contains a
vertex $v$ whose removal disconnects $E$ and separates one triangle
from another.  Let $R$ be the $2$-cell in $D$ corresponding to $v$.
By Lemma~\ref{lem:dual}, $\partial R \cap \partial D$ is disconnected.
Since vertices have degree at most three, the components of $\partial
R \cap \partial D$ are non-trivial paths.
Thus  a portion of $\partial R$  creates a cut path
of length at most $2$.

If $D'$ is connected but singular, at least one of two situations must
occur:  $D'$ contains a vertex of degree~$1$ (a
\emph{spur}) and/or a vertex whose removal disconnects $D'$ (a
\emph{cut vertex}).    Both
possibilities are illustrated in Figure~\ref{fig:double2}.
 
If $v$ is a spur, then there are two edges of $D\setminus D'$ incident
with $v$ that form a cut path separating the $2$-cell $R$ of
$D$ containing these two edges from the rest of the diagram.  Since
$|\partial R| \ge 5$, this creates an exposed path of length $3$
and two consecutive left turns in the turn pattern of $D$.

Vertices in $D'$ have degree at most three,  
and cut vertices have disconnected links, so cut
vertices are always incident with an edge that is not in the
boundary of a $2$-cell of $D'$.  By following 
 a tree of possibilities, we can always find a spur in $D'$ or a vertex $v$
that has a disconnected link and lies in the boundary of a $2$-cell
of $D'$.  There are two $2$-cells in $D \setminus D'$ that
contain $v$, and these $2$-cells must each have $2$ edges incident 
with $2$-cells in the ring $D \setminus D'$.
 Of the two $2$-cells  of $D\setminus D'$
that contain $v$, at least one has $|\partial R \cap D'|
\ge 3$.  Because vertices have degree at most three,
and five edges of $R$ are already accounted for, $R$ must be a hexagon
that shares only one edge with $\partial D$ and the vertices of this
edge yield two consecutive right turns. 
\end{proof}

\begin{cor}\label{cor:double}
If $D$ is a soccer diagram whose turn pattern is $(rl)^i$ and whose
shortest cut path has length at least $3$, then its double dual $D'$
is a soccer diagram.  Moreover, $|\partial D'| = 2i-p$ where $p$ is
the number of pentagons that contain an edge in $\partial D$.
\end{cor}

\begin{proof}
By Lemma~\ref{lem:double}, $D'$ must be both connected and
nonsingular.  Thus it is homeomorphic to a disc and consequently a
soccer diagram.  To prove the second assertion notice that there are
exactly $i$ $2$-cells in the ring of $2$-cells removed from $D$ to
create $D'$.  Moreover, each of these $2$-cells has two of its edges
in $\partial D$ and another two edges are shared with neighboring
$2$-cells in the ring.  Thus each hexagon in the ring contributes two
edges to $\partial D'$ and each pentagon contributes one.
\end{proof}

\begin{cor}[Removing rings]\label{cor:rings}
Let $D$ be a soccer diagram with at most six pentagons.  If the length
of the shortest cut path is $4$ then the double dual of $D$ is a
soccer diagram $D'$ containing exactly six pentagons and $|\partial D|
= |\partial D'|$.
\end{cor}

\begin{proof}
If $D$ contains an exposed path of length $3$ then there would be a
cut path of length less than $4$.  Thus exposed paths of length $3$
cannot exist and by Corollary~\ref{cor:pent}, the turn pattern for $D$
is $(rl)^i$ for some $i$, and $D$ contains exactly six pentagons.  As
a consequence of the left/right alternation, every $2$-cell sharing an
edge with $\partial D$ actual shares two consecutive edges.  Thus
pentagons sharing an edge with $\partial D$ lead to cut paths of
length less than $4$ and all of $2$-cells touching $\partial D$ are
hexagons.  The result now follows immediately from
Corollary~\ref{cor:double}.
\end{proof}

Our key technical result about soccer diagrams is that large diagrams
have long boundary cycles.

\begin{lem}[Large implies long]\label{lem:six}
If $D$ is a soccer diagram with at least six $2$-cells and at most six
pentagons, then $|\partial D| \geq 15$.
\end{lem}

\begin{proof}
We  induct on the number of $2$-cells.  By
Lemma~\ref{lem:min} the statement is true for $k=6$, so suppose it is
true for some $k\geq 6$ and let $D$ be a soccer diagram with exactly
$(k+1)$ $2$-cells.  By Lemma~\ref{lem:cut}, $D$ contains a cut path of
length at most $4$.  Let $P$ be a cut path of minimal length and
consider the number of $2$-cells in $D_1$ and $D_2$.  Without loss of
generality assume that $D_2$ has at least as many $2$-cells as $D_1$.

{\bf Case $1$:} If $D_1$ has at least three $2$-cells, then $D_2$ has
at least four $2$-cells and by Lemma~\ref{lem:min}, $|\partial D_1| +
|\partial D_2| \geq 23$.  Since $|P|\leq 4$ and $|\partial D| =
|\partial D_1| + |\partial D_2| - 2|P|$, we conclude $|\partial D|
\geq 15$.

{\bf Case $2$:} If $D_1$ has two $2$-cells, then $D_2$ has at least
five, and by Lemma~\ref{lem:min}, $|\partial D_1| + |\partial D_2|
\geq 23$ and as in Case~$1$, $|\partial D| \geq 15$.

{\bf Case $3$:} If $D_1$ has one $2$-cell, then $D_2$ has at least
six, and by Lemma~\ref{lem:min}, $|\partial D_1| + |\partial D_2| \geq
20$.  If $|P|\leq 2$, then $|\partial D| \geq 16$.  If $|P|=4$, then
by Corollary~\ref{cor:rings} $\partial D$ has the same length as the
boundary of its double dual $D'$ which is itself a soccer diagram with
at least six $2$-cells.  By induction $|\partial D'| \geq 15$, and the
inequality holds.  Thus we may assume $|P|=3$.  If
$|\partial D_1| > 5$ or $|\partial D_2| > 15$, then $|\partial D| \geq
15$, so we may also assume $D_1$ is a pentagon and $|\partial D_2| =
15$.  This would produce a soccer diagram $D$ with $|\partial D|=14$.
To summarize, the only way for the induction to fail is if there
exists a soccer diagram $D$ with $|\partial D|=14$ whose shortest cut
path has length $3$ and separates a pentagon from the rest of $D$.

If the boundary of this hypothetical diagram contains an exposed
path of length $3$, then the $2$-cell it exposes is a hexagon,
since otherwise there is a cut path of length $2$.  Removing this
hexagon does not change the length of the boundary, and it creates a
diagram $D'$ with $|\partial D'|\geq 15$ by induction.  Thus exposed
paths of length $3$ do not occur.  By Corollary~\ref{cor:pent} the
turn pattern is $(rl)^i$, $D$ contains exactly six pentagons, and
 $i=7$ since $|\partial D|=14$.

Let $D'$ denote the double dual of $D$.  By Lemma~\ref{lem:double},
$D'$ is a soccer diagram with $|\partial D'|\leq 13$ since $D$
contains at least one pentagon which shares an edge with $\partial D$.
By the induction hypothesis $D'$ has at most five $2$-cells and by
Lemma~\ref{lem:min} it actually has at most four.  On the other hand,
the ring removed from $D$ to create $D'$ only contained seven
$2$-cells and at most three of these could be pentagons since the
pentagons must be non-adjacent.  Thus $D'$ contains at least three
pentagons.  Finally, the only soccer diagram with at most four
$2$-cells and at least three pentagons is a hexagon with three
pentagons attached to alternate edges.  Since it is impossible to
reconstruct $D$ by attaching a ring containing a pentagon to this
$D'$, we conclude that $D$ cannot exist.
\end{proof}

Our key result about soccer diagrams follows immediately
from Lemmas~\ref{lem:few} and~\ref{lem:six}. 

\begin{thm}\label{thm:1214}
If $D$ is a soccer diagram with $|\partial D|\leq 14$, at most six
pentagons, and no exposed path of length $4$, then $D$ is one of the
disc diagrams shown in Figure~\ref{fig:1214}.
\end{thm}

%
%
%
%

\section{Algorithm}\label{sec:algorithm}

In this section we review and improve (via a result of Bowditch) the
algorithm for testing the curvature properties of finite piecewise
Euclidean $3$-complexes given in \cite{ElMc-3d}.  For background on
non-positive curvature and piecewise Euclidean complexes, see
\cite{BrHa99}.

A geodesic in a geodesic metric space is \emph{short} if its length is
strictly less than $2\pi$ and \emph{very short} if strictly less than
$\pi$. The original algorithm is based on the standard link condition
for piecewise Euclidean complexes.

\begin{thm}[Link Condition]
\label{thm:linkcond}
A piecewise Euclidean complex is non-positively curved if and only if
the link of each cell has no short closed geodesic.
\end{thm}

Thus deciding whether or not a piecewise Euclidean complex is
non-positively curved depends on checking piecewise spherical
complexes for short geodesics.  Geodesics of this type determine
complexes we call ``circular galleries''.  Rather than give a full
technical definition of a circular gallery, we give a rough definition
and an example.  The reader is referred to \cite{ElMc-alg} and
\cite{ElMc-3d} for precise details.

\begin{defn}[Galleries]
If $\gamma$ is a geodesic in a piecewise spherical complex then the
ordered list of closed simplices through which $\gamma$ passes encodes
a {\em linear gallery determined by $\gamma$}.  If $\gamma$ is a
closed geodesic, this list is given a cyclic rather than a linear
ordering and the result is called the \emph{circular gallery
determined by $\gamma$}.  Linear and circular galleries can also be
determined by paths which are merely close to geodesics.  
\end{defn}

\begin{figure}[htbp]
\begin{tabular}{cc}
\begin{tabular}{c}\includegraphics{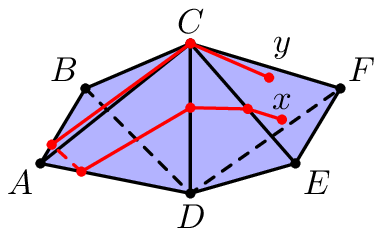}\end{tabular}& 
\begin{tabular}{c}\includegraphics{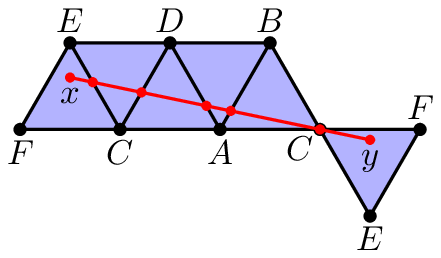}\end{tabular}\\
\end{tabular}
\caption{The $2$-complex and linear gallery described in
Example~\ref{exmp:2-gallery}.
\label{fig:2-gallery}}
\end{figure}

\begin{exmp}\label{exmp:2-gallery}
Let $K$ be the $2$-dimensional piecewise spherical complex formed by
attaching the boundaries of two regular spherical tetrahedra along a
$1$-cell.  The complex $K$ is shown on the left of
Figure~\ref{fig:2-gallery} (where the spherical nature of the
$2$-cells has been left to the reader's imagination).  Let $\gamma$ be
the geodesic shown, which starts at $x$ travels across the front of
$K$, around the back, over the top, and ends at $y$.  The linear
gallery determined by $\gamma$ is shown on the right.
\end{exmp}

In \cite{ElMc-alg} the first two authors prove that given any finite
piecewise Euclidean complex there exists an algorithm to decide
if it is non-positively curved.  In dimension $3$ a second, more
geometric algorithm is available, which has been implemented
as a computer program {\tt cat.g} written in GAP \cite{GAP4}.

The current version of the program is designed to be used with
Euclidean tetrahedra whose edge lengths are square-roots of rationals.
This restriction enables the use of exact arithmetic since all of the
calculations can be carried out in an algebraic number field.  The
program examines the circular galleries that can occur in the the
link of a vertex. The links of other cells in $3$-complexes are easy
to check without a computer.  There are four types of $2$-dimensional
piecewise spherical circular galleries which need to be considered.
If the geodesic passes through a vertex, the gallery is made up of
vertex-to-vertex segments called {\em beads}, which join together to
form a {\em necklace}.  If it does not pass through a vertex, the
gallery it determines is either a disc, an annulus, or a M\"obius band.
Since disc galleries containing short closed geodesics can only exist
in complexes in which the edge links contain short closed geodesics,
these need not be considered.

\begin{figure}[htbp]
\begin{tabular}{ccc}
\begin{tabular}{c}\includegraphics[scale=.75]{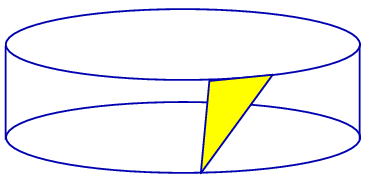}\end{tabular} &
\begin{tabular}{c}\includegraphics[scale=.75]{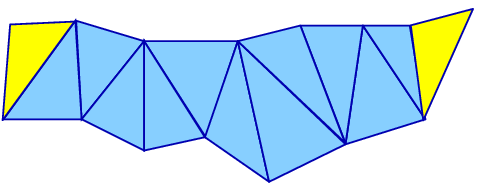}\end{tabular} &
\begin{tabular}{c}\includegraphics[scale=.75]{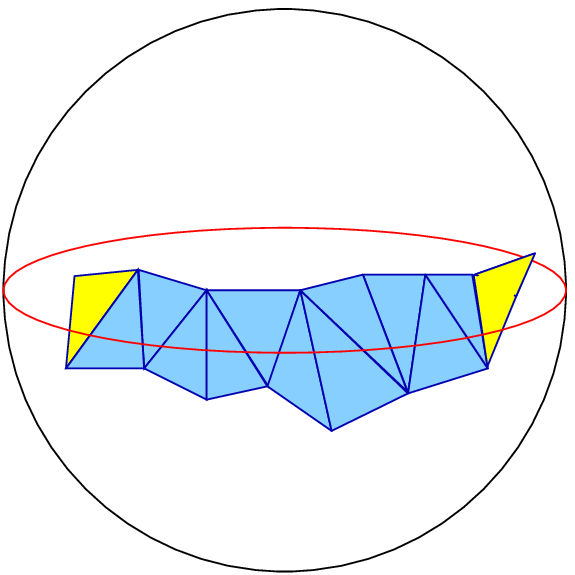}\end{tabular} 
\end{tabular}
\caption{An annular gallery, cut open and developed.}
\label{fig:annular}
\end{figure}

The number of spherical triangles in a circular gallery containing a
short closed geodesic can be bounded ahead of time using only the list
of Euclidean tetrahedra.  Roughly speaking the computer program
proceeds by enumerating every feasible annular or M\"obius gallery and
every bead up to this bound, cuts them open and develops them onto the
$2$-sphere calculating explicit coordinates as it goes.  Then it uses
elementary linear algebra to check each for the existence of a short
closed geodesic.  Details on the algorithm can be found in
\cite{ElMc-3d}; the program is available from the authors'
web-pages.  The cutting open and developing process is illustrated in
Figure~\ref{fig:annular}.

In the remainder of this section we show how the computations
described in \cite{ElMc-3d} can be simplified using a result by Brian
Bowditch \cite{Bo95}. Throughout the remainder of the section let $S$
denote a locally $\cat(1)$ space.  

\begin{defn}[Unshrinkable]
A closed geodesic $\gamma$ in $S$ is  \emph{shrinkable} if there
is a homotopy starting at $\gamma$ such that the length of the closed
curve is non-increasing as a function of time and ends at a curve
whose length is strictly less than its initial length.  Notice that
only the initial curve is required to be a geodesic, so the
equator on a standard metric $2$-sphere is shrinkable.  A closed
geodesic that is not shrinkable is \emph{unshrinkable}.  Even if a
closed geodesic is not shrinkable, there may exist a homotopy such
that the length of the closed curve is unchanging as a function of
time.  In this case we say that the curves at either end are
\emph{equivalent}.
\end{defn}

In this terminology, Bowditch's result 
can be restated as follows:

\begin{lem}[Unshrinkable]\label{lem:shrink}
If $S$ is a locally $\cat(1)$ space that is not globally $\cat(1)$,
then the length of the shortest closed geodesic is the same as the
length of the shortest unshrinkable closed geodesic.
\end{lem}

This leads immediately to the following refinement of the link
condition.

\begin{cor}[Link condition]\label{cor:shrink}
If $S$ is a locally $\cat(1)$ space that does not contain a short
unshrinkable closed geodesic, then $S$ is globally $\cat(1)$.  As a
consequence, a piecewise Euclidean complex is non-positively curved if
and only if the link of each cell has no short unshrinkable closed
geodesic.
\end{cor}

As the next three lemmas show, restricting to unshrinkable
geodesics reduces the number of galleries one needs to inspect.

\begin{lem}[Shrinking annular galleries]\label{lem:annular}
Let $\gamma$ be a closed geodesic in a $2$-dimensional locally
$\cat(1)$ piecewise spherical complex $S$.  If $\gamma$ determines an
annular gallery $\script{G}$, then $\gamma$ is shrinkable.
\end{lem}

\begin{figure}[htbp]
\includegraphics[width=2in]{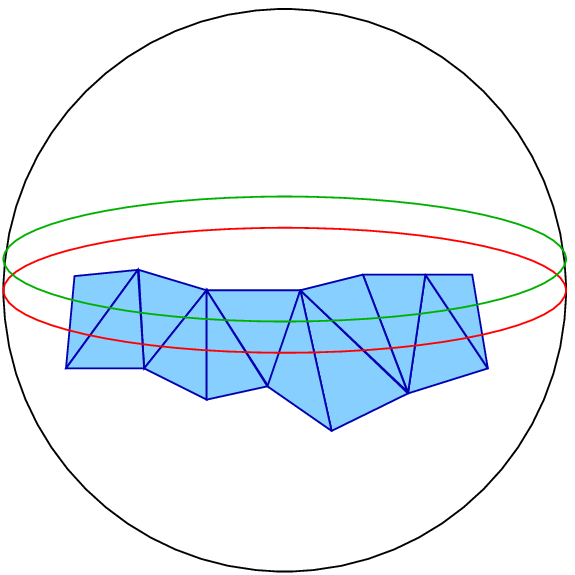}
\hspace*{1cm}
\includegraphics[width=2in]{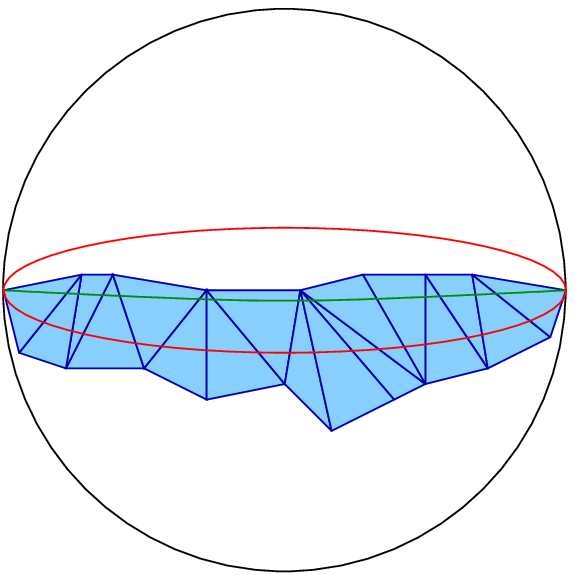}
\caption{Two examples of how geodesics can be tipped}
\label{fig:tip}
\end{figure}

\begin{proof}
When $\script{G}$ is cut open and developed, $\gamma$ is
sent to part of a great circle on the $2$-sphere.  The homotopy which
pushes this path through different lines of latitude shrinks its
length and corresponds in $\script{G}$ to a shorter closed path.
Thus $\gamma$ is shrinkable.  See the left-hand side of
Figure~\ref{fig:tip} for an illustration.
\end{proof}

\begin{lem}[Shrinking M\"obius galleries]\label{lem:moebius}
Let $\gamma$ be a closed geodesic in a $2$-dimensional locally
$\cat(1)$ piecewise spherical complex $S$.  If $\gamma$ determines a
M\"obius gallery and the length of $\gamma$ is at least $\pi$, then
$\gamma$ is shrinkable.
\end{lem}

\begin{proof}
When $\script{G}$ is cut open and developed, $\gamma$ is
sent to path $\gamma'$ in a great circle on the $2$-sphere of length
at least $\pi$.  Let $u$ and $v$ be points in $\gamma'$ that are
antipodal.  The portion of $\gamma'$ between $u$ and $v$ can be
homotoped (in a non-length-changing way) to another geodesic length
$\pi$ connecting them without having the image leave the image of the
cut open M\"obius gallery.  This new path is not locally geodesic at
$u$ or $v$ and can be shortened at either end to produce a strictly
shorter path.  Thus $\gamma$ is shrinkable.  See the right-hand side
of Figure~\ref{fig:tip} for an illustration.
\end{proof}

\begin{lem}[Tipping Beads]\label{lem:vertex}
Let $\gamma$ be a closed geodesic in a $2$-dimensional locally
$\cat(1)$ piecewise spherical complex $S$.  If $\gamma$ determines a
necklace gallery $\script{G}$ such that a single bead contains a
portion of $\gamma$ of length strictly more than $\pi$, then $\gamma$
is shrinkable.  Moreover, if $\gamma$ determines a necklace gallery
$\script{G}$ that contains a single bead containing a portion of
$\gamma$ of length exactly $\pi$, then $\gamma$ is equivalent to a
path $\gamma'$ which determines a necklace gallery in which all beads
contain strictly less than $\pi$ of $\gamma'$.
\end{lem}

\begin{proof}
If there is a bead containing more than $\pi$ of $\gamma$, then we can
pick $u$ and $v$ in the interior of the bead which are connected by a
portion of $\gamma$ of length exactly $\pi$.  The rest of the proof
mimics the proof of Lemma~\ref{lem:moebius}.  If there is a bead
containing a portion of $\gamma$ of length exactly $\pi$, then we pick
the vertices at either end through which $\gamma$ passes as our $u$ and
$v$ and proceed to tip the portion of $\gamma$ between them.  Because
the link of $u$ may have excess curvature, we do not know that
$\gamma$ can be locally shortened after $\gamma$ has been tipped.  On
the other hand, we can continue tipping the portion of $\gamma$
between $u$ and $v$ until the path hits the boundary cycle of the
bead.  Because this boundary cycle is a piecewise geodesic path, the
portion of its boundary included in the tipped path must include a new
vertex of $\script{G}$.  Since $\gamma$ is assumed to be short, there
is only one bead of length $\pi$, and in this equivalent path this
one long bead has been broken up into at least two shorter ones.
\end{proof}

The original program searched for all annular and M\"obius galleries
that contain short closed geodesics and all beads which contain
geodesics of length less than $2\pi$, which are then strung together to
form necklaces.  By the last three lemmas we do not need to search for
annular galleries at all, or for the longer types of beads and
M\"obius galleries.  To appreciate the magnitude of this
simplification see Remark~\ref{rem:runtime}.

\section{Metric and output}\label{sec:output}

We begin this section by defining the shapes used
to give $5/6^*$-triangulated $3$-manifolds 
piecewise Euclidean structures.

\begin{defn}[The metric]\label{def:metric}
Let $M$ be a $5/6^*$-triangulated, closed $3$-manifold.  We 
make $M$ a \emph{metric $5/6^*$-triangulated $3$-manifold}
by assigning a 
length of $\sqrt{3}$ to each edge of degree $6$, a length of $2$ to
each edge of degree $5$, and  metrics of
 the unique  Euclidean simplices whose edge
lengths match those assigned to their $1$-skeletons to the
triangles and tetrahedra.  
\end{defn}

\begin{figure}[htbp]
\includegraphics[scale=1.2]{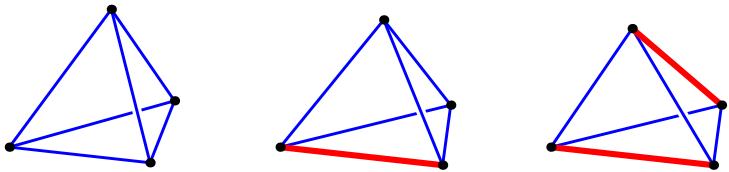}
\caption{The three types of tetrahedra (regular, mixed, and Coxeter)
in a $5/6^*$-triangulation.  The thin edges have length~$\sqrt{3}$ and
the thick edges have length~$2$.\label{fig:shapes}}
\end{figure}

\begin{defn}[The tetrahedra]
The fact that edges of degree $5$ cannot belong to the same $2$-cell
means that there are only three equivalence classes of metric
tetrahedra in $M$: those with $0,1$ or $2$ edges of degree $5$.  See
Figure~\ref{fig:shapes}.  Notice that the first tetrahedron is regular
and the third one is a Coxeter shape with dihedral angles 
$\pi/2$ and $\pi/3$ around the edges of degree $5$ and $6$, respectively.
We refer to these three metric tetrahedra as \emph{regular},
\emph{mixed} and \emph{Coxeter} tetrahedra.
\end{defn}

Straightforward computations show:

\begin{lem}[Dihedral angles]\label{lem:dihedral}
Let $M$ be a closed $3$-manifold with a metric $5/6^*$-triangulation
and let $e$ be an edge in a tetrahedra $T$ in $M$.  If $e$ has degree
$5$, then the dihedral angle in $T$ at $e$ is more than $2\pi/5$ and
at most $\pi/2$.  If $e$ has degree $6$, then the dihedral angle in
$T$ at $e$ is at least $\pi/3$, strictly less than $\pi/2$, and equal
to $\pi/3$ only if $T$ is Coxeter. As a consequence, the links of
edges in $M$ are metric circles of length at least $2\pi$ and exactly
$2\pi$ if and only if $e$ has degree $6$ and is surrounded by six
Coxeter tetrahedra.
\end{lem}

The software {\tt cat.g} returns a list of $4$ trivial beads, $71$
non-trivial beads and $144$ M\"obius galleries, for these three
tetrahedra.  Because we are interested in triangulated $3$-manifolds,
and M\"obius strips cannot be immersed into $2$-spheres, the M\"obius
strips in the output can be safely ignored.

\begin{figure}[htbp]
\includegraphics[scale=1.4]{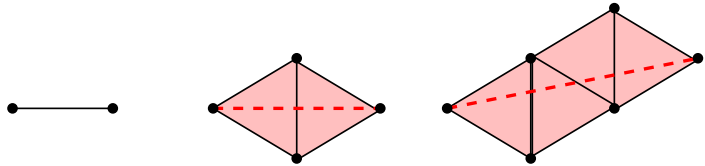}
\caption{The three types of beads.}
\label{fig:beads}
\end{figure}

\begin{defn}[Bead types]\label{def:beads}
Combinatorially all $75$ metric beads of length less than $\pi$ in the
output look like the one of the three non-metric beads shown in
Figure~\ref{fig:beads}.  The differences come from the metrics.
Specifically there are $4$ different metric edges, $26$ metric beads
consisting of two triangles, and $45$ metric beads consisting of four
triangles.  We refer to these underlying combinatorial structures
as beads of \emph{type~$A$}, \emph{type~$B$} and \emph{type~$C$}
respectively.  In each case, the geodesic contained in the bead starts
at the leftmost vertex and ends at the rightmost vertex.  

Similarly, the \emph{type of a necklace}, refers to its underlying
combinatorial structure, rather than its metric.  A necklace
containing at least one non-trivial bead is called a \emph{thick
necklace} while those consisting solely of trivial beads are called
\emph{thin necklaces}.
\end{defn}

The following estimates on the lengths of
geodesics in each of the three bead types
immediately implies Corollary~\ref{cor:necklaces}.

\begin{lem}[Lower bounds]\label{lem:lower-bounds}
The length of a geodesic in a bead of type $A$, $B$ or $C$ is bounded
below by $.304\pi$, $.5\pi$, and $.832\pi$, respectively.
\end{lem}

\begin{cor}[Necklaces]\label{cor:necklaces}
There are $28$ types of necklaces consisting of short beads that might
contain a closed geodesic of length less than $2\pi$.  These $28$
possibilities are
$A$,      $A^2$,     $A^3$,  $A^4$,    $A^5$,
$A^6$,    $B$,       $AB$,   $A^2B$,   $A^3B$, 
$A^4B$,   $B^2$,     $AB^2$, $A^2B^2$, $ABAB$, 
$A^3B^2$, $A^2BAB$,  $C$,    $AC$,     $A^2C$, 
$A^3C$,   $BC$,      $ABC$,  $A^2BC$,  $ABAC$, 
$B^2C$,   $C^2$, and $AC^2$.
\end{cor}

\begin{rem}[Speed]\label{rem:runtime}
The entire computation took less than one hour to complete on a 850MHz
PC running GAP under Linux.  It tested approximately $110$ thousand
galleries and yielded the $75$ beads described above.  As an
indication of the benefits of the simplification described in
Section~\ref{sec:algorithm}, we note that an earlier computer search
--- using the original algorithm --- took more than $2$ months to
complete, tested approximately $300$ million galleries and yielded $
\simeq 12\,000$ beads.  Moreover, instead of the three combinatorial
types shown in Figure~\ref{fig:beads}, there were more than $100$
combinatorial types.  In other words, the restriction to unshrinkable
geodesics not only improved the length of the computation, it also
greatly simplified the output.
\end{rem}

\section{Main result}\label{sec:main}
Throughout this section let $M$ be a metric $5/6^*$-triangulation of a
$3$-manifold, let $S$ be the link of a $0$-cell in $M$, and let
$\gamma$ be a closed geodesic contained in $S$ which determines a
necklace gallery $\script{G}$ consisting of short beads.

\begin{figure}
\includegraphics[scale=1.4]{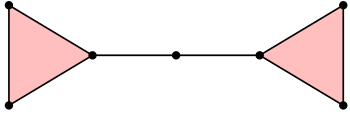}
\caption{A barbell with $3$ cut vertices.\label{fig:barbell}}
\end{figure}

\begin{defn}[Barbells]\label{def:barbells}
Let $\script{G}$ be a thick necklace gallery.  The linear gallery
consisting of a (possibly empty) sequence of trivial beads with a
single triangle on either end is called a \emph{barbell}.  See
Figure~\ref{fig:barbell}.  Notice that a barbell containing $i-1$
trivial beads contains exactly $i$ cut vertices.  Also note that the
transition between two non-trivial beads leads to a barbell consisting
of two triangles joined at a vertex.
\end{defn}

\begin{defn}[Good perturbations]\label{def:perturb}
If $\gamma$ is perturbed slightly so that it avoids all of the
vertices in $S$, then the circular gallery determined by the new path
is an annular gallery.  In order to maintain control over the
combinatorial properties of the annular gallery which results, we
define a \emph{good perturbation} of $\gamma$ as follows.  If
$\script{G}$ is a thin necklace, then $\gamma'$ is the boundary
curve of an $\epsilon$-neighborhood of $\gamma$ in $S$ that passes
through the minimum number of triangles.  If both boundary curves pass
through the same number of triangles, both boundary curves are
good perturbations.  If $\script{G}$ is a thick necklace, then 
define $\gamma'$ one barbell at a time.  In each barbell we
consistently push the path $\gamma'$ to the left of all the cut
vertices or to the right of all the cut vertices, whichever minimizes
the number of triangles through which it passes.  See
Example~\ref{exmp:resolve}.
\end{defn}

\begin{figure}[htbp]
\includegraphics[scale=1.4]{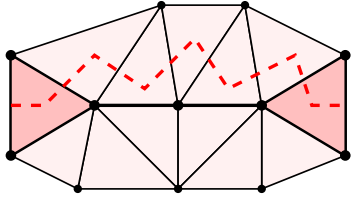}
\caption{A barbell, its neighborhood and a perturbed geodesic.}
\label{fig:barbell2}
\end{figure}

\begin{exmp}\label{exmp:resolve}
Figure~\ref{fig:barbell2} shows a barbell in which the good
perturbation lies above the three cut vertices.  Since we are only
interested in the gallery determined by the perturbed path (and not
the path itself), the jagged line connecting the centers of the
triangles is a perfectly good representative of the perturbed path.
\end{exmp}

We need two technical lemmas about the paths determined by good
perturbations; one for thin necklaces and one for thick necklaces.

\begin{lem}[Thin necklaces]\label{lem:thin}
Let $\gamma$ be a geodesic of length at most $2\pi$ which determines a
thin necklace, let $\gamma'$ be a good perturbation of $\gamma$ and
let $P$ be the closed immersed path that is the internal dual of the
annular gallery determined by $\gamma'$.  If $c$ is the absolute value
of the combinatorial turning angle for $P$, then $|P| + c \leq 12$.
Moreover, the turn pattern for $P$ does not contain three consecutive
left turns or three consecutive right turns.
\end{lem}

\begin{proof}
Since $\script{G}$ is thin it has type $A^i$ for some $i$ and by
Lemma~\ref{lem:lower-bounds}, $i\leq 6$.  Notice that $i$ is also the
number of vertices in $\gamma$.  For each vertex $v$ in $\gamma$, the
two boundary curves of an $\epsilon$-neighborhood of $\gamma$ pass
through five or six corners of triangles.  The triangles on either
side of the trivial beads of $\gamma$ are double counted in the sense
that the boundary curves traverse two corners of these triangles in a
row before moving on to a new triangle.  They cannot traverse three
corners in a row since this would require $\gamma$ to traverse two
sides of the triangle consecutively, and this is prohibited by the
size of the dihedral angles (Lemma~\ref{lem:dihedral}) and the fact
that $\gamma$ is a local geodesic.  Thus the two paths together
determine annular galleries that contain at most $4i$ triangles
combined, and one annular gallery has at most $2i$ triangles.

Without loss of generality assume that the good perturbation is the
one to the left of $\gamma$ as $\gamma$ is traversed.  Let $R$ be a
triangle that contains a portion of $\gamma'$ and notice that $R$
corresponds to a vertex in $P$ corresponding to a left turn if and
only if $\partial R$ contains an edge of $\gamma$.  This implies that
the turn pattern for $P$ contains exactly $i$ left turns and that $|P|
+ c = (n_l+n_r) + (n_l-n_r) = 2n_l = 2i \leq 12$.

Suppose $P$ contains three consecutive right turns.  The three
triangles corresponding to these turns, plus the triangle immediately
before and after (as traced out by $\gamma'$), all contain a common
vertex $v$ in their boundaries.  Since $v$ has degree at most $6$, the
path $\gamma$ would have to make a sharp turn at $v$ and $\gamma$
would not be a local geodesic.  If $P$ contains three consecutive left
turns, then there are three consecutive edges in $\gamma$ (with
separating vertices $u$ and $v$) such that all three triangles to the
left of these edges contain a common apex.  In order for $\gamma$ to
be a local geodesic at $u$ and $v$, by Lemma~\ref{lem:dihedral} both
vertices would need to have degree~$5$, which is forbidden.
\end{proof}

\begin{lem}[Thick necklaces]\label{lem:thick}
Let $\gamma$ be a short unshrinkable geodesic which determines a thick
necklace consisting of short beads, let $\gamma'$ be a good
perturbation of $\gamma$, and let $P$ be the closed immersed path that
is the internal dual of the annular gallery determined by $\gamma'$.
If $c$ is the absolute value of the combinatorial turning angle for
$P$, then $|P| + c \leq 14$.  Moreover, the turn pattern for $P$ does
not contain three consecutive left turns or three consecutive right
turns.
\end{lem}

\begin{proof}
The proof is similar to Lemma~\ref{lem:thin}, but it proceeds one
barbell at a time.  Consider a barbell in $\script{G}$ with $i$ cut
vertices.  For each cut vertex $v$ in the barbell, the two possible
perturbations of $\gamma$ pass through five or six corners of
triangles.  As before the triangles on either side of the trivial
beads are double counted in the sense that the boundary curves
traverse two corners of these triangles in a row before moving on to a
new triangle.  They cannot traverse three corners in a row since this
would require $\gamma$ to traverse two sides of the triangle
consecutively, which is prohibited by the size of the dihedral angles
(Lemma~\ref{lem:dihedral}) and the fact that $\gamma$ is a local
geodesic.  Since the triangles at either end already existed, the two
possible perturbations pass through at most $6i-2(i-1)-2 = 4i$ new
triangles.  In particular, one of them passes through at most $2i$ new
triangles.

Without loss of generality assume that the good perturbation is the
one to the left of the portion of $\gamma$ in this barbell as $\gamma$
is traversed.  Let $R$ be a triangle which is traversed by this
portion of $\gamma'$ and notice that $R$ corresponds to a vertex in
$P$ which is a left turn if and only if $\partial R$ contains an
trivial bead of $\gamma$ or $R$ is one of the original two triangles
in the barbell.  This implies that this portion of the turn pattern
for $P$ contains exactly $i$ left turns.  As above, the number of
vertices in this portion of $P$ plus the absolute value of the
combinatorial turning angle for this portion is twice the number of
left turns for this portion, which is $2i$.  Finally, notice that the
only portions of $P$ that are not contained in some barbell are the
portions corresponding to the two interior triangles in a bead of
type~$C$.  Since these two triangles become vertices in $P$ which turn
in opposite directions, they contribute to the length of $P$ but not
to $c$.  It is now routine to calculate that $14$ is an upper bound
for $|P|+c$ for each of the $22$ thick cases listed in
Corollary~\ref{cor:necklaces}.  Finally, the arguments that $P$ has no
three consecutive left [right] turns is identical to the one given
above and is omitted.
\end{proof}

It is now relatively easy to show that the paths that good
perturbations determine are in fact simple.

\begin{lem}[Simple]\label{lem:simple}
Let $\gamma$ be either a geodesic of length at most $2\pi$ which
determines a thin necklace or a short geodesic which determines a
thick necklace consisting of short beads.  If $P$ is the closed
immersed path which is the internal dual of the annular gallery
determined by a good perturbation of $\gamma$, then $P$ is a simple
closed path.
\end{lem}

\begin{proof}
Suppose $P$ is a not embedded and let $Q$ be a closed subpath of $P$
of minimal length.  Since $|P|\leq 14$ by Lemma~\ref{lem:thin} and
Lemma~\ref{lem:thick}, $|Q|\leq 7$.  Since $Q$ itself is embedded it
divides the soccer tiling into two soccer diagrams $D_1$ and $D_2$,
one of which, say $D_1$, has fewer than six pentagons.  By
Lemma~\ref{lem:min}, $D_1$ consists of a single $2$-cell and thus $P$
contains at least four consecutive left [right] turns, contradicting
Lemma~\ref{lem:thin} or Lemma~\ref{lem:thick}.
\end{proof}

We can now show that our hypothetical short closed unshrinkable
geodesic $\gamma$ does not exist.

\begin{lem}[Vertex links]\label{lem:vertex-links}  
If $S$ is the link of a vertex in a metric $5/6^*$-triangulated,
closed $3$-manifold, then it does not contain any short closed
unshrinkable geodesics.  In addition, $S$ does not contain any closed
geodesics of length $2\pi$ in its $1$-skeleton.
\end{lem}

\begin{proof}
If $S$ contains a short closed unshrinkable geodesic, then by
Lemma~\ref{lem:vertex} it also contains a short closed unshrinkable
geodesic which determines a necklace consisting of short beads.  Let
$\gamma$ be either a closed geodesic in $S$ of length at most $2\pi$
which determines a thin necklace or a short closed geodesic in $S$
which determines a thick necklace consisting of short beads.  By
Lemmas~\ref{lem:thin},~\ref{lem:thick}, and~\ref{lem:simple} there is
a perturbation of $\gamma$ that determines an annular gallery whose
internal dual is a simple closed path $P$ in the dual soccer tiling.
Moreover, if $c$ denotes the absolute value of the combinatorial
turning angle for $P$, then $|P|+c\leq 14$ and the turn pattern for
$P$ does not contain three consecutive left turns or three consecutive
right turns.  Since $P$ is embedded it divides the soccer tiling into
two soccer diagrams $D_1$ and $D_2$, one of which, say $D_1$, has
fewer than six pentagons.  Since $P$, with the appropriate orientation
is $\partial D_1$, Theorem~\ref{thm:1214} implies that $D_1$ is one of
the two soccer diagrams in Figure~\ref{fig:1214}.  Since both diagrams
have $c=4$ and $|\partial D| + c > 14$, we have a contradiction.
\end{proof}

The final property that we need to establish is the following.

\begin{lem}[No flat planes]\label{lem:no-flats}
If $M$ is a closed $3$-manifold with a metric $5/6^*$-triangulation,
then the universal cover $\widetilde{M}$ does not contain any
isometrically embedded flat planes.
\end{lem}

\begin{figure}
\bt{c}\includegraphics[width=3cm,height=3cm]{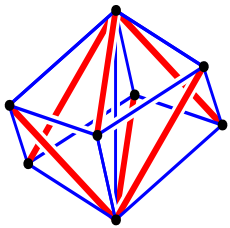}\et
\caption{A cage around a degree $6$ edge.}
\label{fig:cage}
\end{figure}

\begin{proof}
Let $\phi:\R^2\to \widetilde M$ be an isometric embedding of a flat
plane and let $F=\phi(\R^2)$.  If $F$ is transverse to an edge $e$ in
$\widetilde M$ and $x$ is the unique point in $e \cap F$, then the
link of $x$ in $F$ is a closed geodesic loop $\gamma$ of length $2\pi$
in the space of directions of $x$ in $\widetilde M$.  Since the space
of directions of $x$ is an orthogonal join of $\sph^0$ with the metric
circle which is the link of $e$ in $\widetilde M$, the space of
directions of $x$ is either a standard $2$-sphere (when the link of
$e$ has length exactly $2\pi$) or a branched cover of $\sph^2$ around
antipodal points (when the link of $e$ has length greater than
$2\pi$).  Notice that in either case the loop $\gamma$ must avoid the
points corresponding to the edge $e$ since $F$ is transverse to $e$.
In the branched case this is impossible and in the non-branched case,
we can assume (Lemma~\ref{lem:dihedral}) that $e$ has degree $6$ and
is surrounded by six Coxeter tetrahedra.  This shows that $F$ cannot
cross any edge of degree~$5$.  Figure~\ref{fig:cage} shows an edge $e$
of degree $6$ surrounded by six Coxeter tetrahedra.  Regardless of the
location of the point $x$ in $e$, there does not exist a portion of a
flat plane through $x$ which does not extend through one of the six
edges of degree $5$.  The key observation is that there are three
edges of degree $5$ extending down from the top of $e$, another three
edges of degree $5$ extending up from the bottom of $e$, and these six
edges interleave.  In other words, the edges of degree $5$ form a cage
from which a portion of a flat plane cannot escape.

The remaining possibility is that $F$ is completely contained in the
$2$-skeleton of $\widetilde M$.  If $v$ is a vertex of $\widetilde M$
contained in $F$, then the link of $v$ in $F$ corresponds to a closed
geodesic loop of length $2\pi$ in the $1$-skeleton of the link of $v$
in $\widetilde M$.  By Lemma~\ref{lem:vertex-links} this is also
impossible.
\end{proof}

\renewcommand{\thethm}{\ref{thm:main}}
\begin{thm}[Main Theorem]
Every $5/6^*$-triangulation of a closed $3$-manifold $M$ admits a
piecewise Euclidean metric of non-positive curvature, where the
universal cover $\widetilde{M}$ contains no isometrically embedded
flat planes.  As a consequence, $\pi_1(M)$ is word-hyperbolic.
\end{thm}

\begin{proof}
Let $M$ be such a piecewise Euclidean $3$-manifold.  By
Theorem~\ref{thm:linkcond} it is sufficient to show that the links of
the cells in $M$ do not contain short geodesic loops.  This is
trivially true for the links of $3$-cells and $2$-cells (whose links
are empty and discrete), and it is also true for links of $1$-cells
and $0$-cells by Lemma~\ref{lem:dihedral} and
Lemma~\ref{lem:vertex-links}.  Thus $M$ is non-positively curved.
Since there are no isometrically embedded flat planes in
$\widetilde{M}$ by Lemma~\ref{lem:no-flats}, the final assertion
follows immediately from Theorem~III.$\Gamma$.3.1 in \cite{BrHa99}.
\end{proof}

\def\cprime{$'$}

\end{document}